\documentclass[preprint,12pt]{article}

\usepackage{amssymb}
\usepackage{graphicx}

\begin{document}

\vspace{-1cm}

\title{Circadian rhythm and cell population growth}
\maketitle
\author{Jean Clairambault$^{\footnotemark[1],\footnotemark[2]}$,
  St\'ephane Gaubert$^{\footnotemark[3],\footnotemark[4]}$ and Thomas Lepoutre$^{\footnotemark[1],\footnotemark[5]}$}
\newtheorem{Theorem}{Theorem}
\newtheorem{lemma}{Lemma}[Theorem]
\newtheorem{cor}{Corollary}[Theorem]


\footnotetext[1]{ INRIA, projet BANG, Domaine de Voluceau, BP 105, 78156 Le Chesnay Cedex France}
\footnotetext[2]{INSERM U 776,
H\^opital Paul-Brousse,
14, Av. Paul-Vaillant-Couturier
F94807 Villejuif cedex}

\footnotetext[3]{ INRIA Saclay -- \^Ile-de-France, projet MAXPLUS} 
\footnotetext[4]{ CMAP, Ecole Polytechnique, 91128 Palaiseau Cedex, France}
\footnotetext[5]{ UPMC Univ Paris 06, UMR 7598, Laboratoire Jacques-Louis Lions, F-75005, Paris, France.}
\begin{abstract}
\footnotesize{
Molecular circadian clocks, that are found in all nucleated cells of mammals, are known to dictate rhythms of approximately 24 hours ({\em circa diem}) to many physiological processes. This includes metabolism (e.g., temperature, hormonal blood levels) and cell proliferation. It has been observed in tumor-bearing laboratory rodents that a severe disruption of these physiological rhythms results in accelerated tumor growth. 

The question of accurately representing  the control exerted by circadian clocks on healthy and tumour tissue proliferation to explain this phenomenon has given rise to mathematical developments, which we review. The main goal of these previous works was to examine the influence of a periodic control on the cell division cycle in physiologically structured cell populations, comparing the effects of periodic control with no control, and of different periodic controls between them. We state here a general convexity result that may  give a theoretical justification to the concept of cancer chronotherapeutics. Our result also leads us to hypothesize that the above mentioned effect of disruption of circadian rhythms on tumor growth enhancement  is indirect, that, is this enhancement is likely to result from the weakening of healthy tissue that are at work fighting tumor growth.\\

\textbf{keywords :}
mathematical biology, partial differential equations, physiologically structured models, growth processes, eigenvalues , circadian rhythms, cancer}
\end{abstract}


\section{Introduction: a challenging question from biology}
\label{intro}

The existence of circadian rhythms in humans has been known for centuries \cite{lemmer}, but only recently, in the last thirty years, has their molecular nature been located and understood in cell physiological mechanisms \cite{reppert-weaver2002, hastings2003}. Circadian clocks (from the latin {\em circa diem}, about one day) have been shown to be present in all nucleated cells, and to be conducted by a central circadian clock. This clock consists in about 20000 coupled neurons located in the suprachiasmatic nuclei of the mammalian hypothalamus, and is itself reset by external light through the retinohypothalamic tract. Circadian clocks influence by nervous or hormonal messengers cell metabolism and tissue proliferation \cite{reppert-weaver2002, bjarnason1999, bjarnason2002, fu-lee2003, levi-schibler2007}.

In a series of papers \cite{filipski2002, filipski2004, filipski2005} reporting results from biological experiments on laboratory rodents, Filipski et al. have shown that a severe disruption, obtained either by surgery or by light entrainment perturbations, of the central hypothalamic circadian clock in tumor-bearing mice leads to an acceleration of tumor growth. These experiments were conducted to give experimental confirmation of the fact, known in the clinic of cancer, that patients with maintained circadian rhythms, i.e., significant amplitude of 24h-periodic signals of rest-activity rhythm, temperature, and blood cortisol level, have a much better prognosis than those whose circadian rhythms are damped or ablated  \cite{levi-schibler2007}.

From a mathematical point of view, this question is related to the control of growth processes. As far as linear models (or locally linearized models of more complex ones) are concerned, the natural output of such growth processes to be controlled is a dominant eigenvalue. Tissue proliferation observed at the macroscopic level, e.g., of tumor growth, relies at the microscopic level on cell division in populations of cells. For this reason,the first author, and others, in 2003 initiated a series of studies    \cite{ clmp03, CMP, CMP2, CGP, Bekkal1, Bekkal2, Doumic, Lepoutre_MMNP}. In these studies,  they designed  physiologically based models of the cell division cycle in proliferating cell populations, and of their external control. In large cell populations, it is natural to choose for these models systems of partial differential equations including parameters that will become targets for an external control. Such control may be physiological, i.e., hormonal or circadian, or else pharmacological, by drugs acting on the cell division cycle, as is the case in the clinic of cancer (cytotoxic drugs). It is noteworthy that {\em chronotherapeutics of cancer} precisely uses with clinical success periodic pharmacological control of the cell division cycle \cite{levi-schibler2007, levi2001}. But even without considering pharmacological control, and only physiological control by circadian inputs, we have been led to examine the effects of a periodic control on the cell division cycle in proliferating cell populations.

Underlying biological questions on which mathematical modelling can shed some light are: What is the exact effect of a periodic control on cell proliferation measured by a dominant eigenvalue of the process?  In what sense  do circadian rhythms act on cell and tissue proliferation? Enhancement or lowering? Is the aforementioned biological phenomenon observed in tumor-bearing mice the result of an effect on tumor cells, or on healthy cells hampered in their fight against cancer?  Can an enhancement of circadian (hormonal, photic)  rhythms by artificial external delivery be used as an adjuvant treatment against cancer? Mathematical models can certainly not answer all these questions, but they may give guidelines to help solve them, and this has been for us an incentive to undertake studies on periodic control of the cell division cycle.

In this paper, we review our results, including a new unifying convexity inequality for the Floquet eigenvalue (Theorem~\ref{theorem:convexityth}). The latter inequality implies that, assuming the therapeutic control only influences the death rates of cells, the stationary control can always be replaced by a periodic control with the same average, in such a way that the growth rate is increased. This property accounts for the importance of circadian effects on toxicity for healthy tissue. The fact that the growth rate is generically increased by a periodic control should be compared with the experimental observation of Filipski et al. (\cite{filipski2004,filipski2005}) that tumor proliferation is decreased by a periodic control. In this perspective, our results support the hypothesis that the effect of circadian control on tumor proliferation is likely to be indirect, not resulting as firstly envisioned from a direct action on tumor cells, but more probably from the weakening of healthy tissues that are at work fighting tumor growth.

This paper is organized as follows : in section 2 we detail the foundations of our type of models. In section 3, we review our results on the dominant eigenvalue, state our main new convexity result and explain how these results are supporting cancer chronotherapeutics and lead us to change our point of view on  the initial observation of Filipski et al. on tumor growth enhancement by circadian disruption as explained in the last paragraph. 
\section{Physiological and pharmacological control on cell proliferation}
\label{circclockscellprolif}

The cell division cycle can be modelled in proliferating cell populations by age-structured partial differential equations. Note that in this perspective, space is not necessarily a relevant variable. One can reasonably assume that tissue vascularization, both in healthy tissues and in evolved tumors, as one can observe under a microscope, is extended enough so as to allow us to hypothesize a homogeneous distribution of drugs and hormones from the central blood compartment into a common cell population, healthy or tumoral. We emphasize here that the models we use claim not to reflect the position of a plain observer of solid tumor growth, tuning biophysical parameters, but rather the (supposedly more active) role of a pharmacologist acting by known molecular inputs - drugs or hormones, or other control mechanisms - on physiological targets in a reduced model of cell proliferation, with the aim to limit cancer proliferation. It is having this ``physician rather than physicist'' caveat in mind that we have designed such physiologically based models where age in cell cycle phases, not space, is the main structure variable.

The cell division cycle is classically divided in four phases, namely $G_1$, $S$ (for DNA synthesis), $G_2$ and $M$ (for mitosis, i.e., effective cell separation into two daughter cells).  In each of these phases, that may be considered as subpopulations of the total considered cell population, cells proceed unidirectionally toward the next phase, with cell doubling at mitosis, along an age axis. This physiological variable age in fact represents different biological processes such as DNA synthesis for $S$ phase, and proteic synthesis for the growth phases $G_1$ and $G_2$ during which cells prepare the biological material necessary for DNA duplication (in $S$ phase) or microtubule synthesis and assembling (in $M$ phase). In each of these phases, progression speed (i.e., age versus sideral time), death rate within the phase and rate of transition to the next phase are physiological control targets. Progression speed may be enhanced by external growth factors that bind to membrane receptors. Death rate is controlled by pro- and antiapoptotic factors (apoptosis being programmed  cell death, necessary to maintain physiological equilibria incell populations). And transition rates are controlled by proteic complexes of cyclins and cyclin dependent kinases that show stiff dynamics and may be thought of as gates, allowing or not passage from one phase to the next.

If one does not take in the present studies growth factors into account, we can  simplify these physiological settings by setting progression speed in the cell division cycle to a constant.  The controlled parameters of the population dynamics partial differential equation (PDE) model described below are then death rates, hereafter $d_i(t,x)$, and transition rates, hereafter $K_{i\rightarrow i+1}(t,x)$, in phase $i, (1\leq i\leq I)$. In a first study on the dependency of the first eigenvalue of the growth process as a function of periodic control on death rates or on transition rates \cite{CMP, CMP2}, had been shown a result that was rather a surprise to the authors. Biological evidence was indeed suggesting the exact opposite of what was then obtained: time-periodic  death rates, as opposed to non controlled, i.e., constant death rates, with the same arithmetic means, yield higher first eigenvalues, and thus enhanced proliferation of the population (detailed below). Assuming that periodic control represents the normal circadian rhythm regularly reset by dark-light alternance, and no control (i.e., constant death rates) a disrupted circadian clock, we were expecting from the aforementioned experiments on laboratory rodents the converse, i.e., higher first eigenvalues in the case of a disrupted clock. This apparent discrepancy was the motivation for us to understand more precisely the dependency of eigenvalues on a periodic control on cell division processes. We briefly review below the previously obtained results, and generalize them by giving a new result relying on a convexity argument.

Moreover, we consider in a simple setting the representation of chronotherapeutics by a periodic effect on death rates and we show that  chronotherapeutics is bound to yield generically better results than constant drug delivery.

\section{Periodic control on the cell division cycle: modelling and results}
\subsection{The linear age-structured model and the dominant eigenvalue}
\label{periodiccontrol}
There are many ways to model the circadian control on the cell cycle. One of the simplest is to consider equations with time-periodic coefficients. 
We have chosen to consider the framework of renewal equations. To study the effect of periodic forcing, we have studied its effects on linear models. The typical model, introduced in \cite{clmp03} is built as follows:
\begin{itemize}
\item the cell cycle is divided into $I$ successive physiological phases (typically $I=4$ and the phases are as usual $G_1, S, G_2, M$), 
\item in each phase $i$ the population number (or density) is represented by a cell variable $n_i(t,x)$ ($t$ is time, $x$ is the age in the phase),
\item cells can leave each phase $i$ with a rate $K_{i\rightarrow i+1}$ to the next phase $i+1$,  which they enter with age $0$. 
\end{itemize}
This can be expressed as
\begin{equation}\label{eq:cycle_cellulaire}
\left\lbrace\begin{array}{l}
\partial_t n_i(t,x)+\partial_x n_i(t,x)+\left[d_i(t,x)+K_{i\rightarrow i+1}(t,x)\right]n_i(t,x)=0,\\[0.3cm]
n_{i+1}(t,0)=\int_0^\infty K_{i\rightarrow i+1}(t,x)n_i(t,x)dx,\\[0.3cm]
n_1(t,0)=2\int_0^\infty K_{I\rightarrow 1}(t,x)n_I(t,x)dx.\end{array}\right.
\end{equation}
Note that every coefficient is taken as time-dependent. More precisely, we assume coefficients to be $T-$periodic with respect to the time variable $t$. Such systems are characterized by a dominant eigenvalue $\lambda_F$ ($F$ for Floquet, referring to Floquet theory for the analysis of periodic solutions of differentiable systems, see appendix (A.3) for a precise definition) which governs the growth behavior of solutions, in the sense that solutions may be expressed as  $e^{\lambda_F t}$ times a bounded term. In the simplest case $I=1$, without death rates, already many mathematical effects of time heterogeneity can be observed. The equation then reads
\begin{equation}\label{eq:division}\left\lbrace\begin{array}{l}
\partial_t n(t,x)+\partial_x n(t,x)+K(t,x)n(t,x)=0,\\[0.3cm]
n(t,0)=2\int_0^\infty K(t,x)n(t,x)dx.\end{array}\right.
\end{equation}
Note that an asymptotic link with a discrete system can be established  and that in such discrete systems, in some situations, a paradoxical decrease of the growth rate may be obtained by increasing the division rate $K(t,x)$ \cite{Lepoutre_PHD}.
For a more complete description of the dominant eigenvalue and its associated eigenvectors, we refer to \cite{CMP,CGP,Lepoutre_MMNP,Lepoutre_PHD,MMP}.

\subsection{Averaged coefficients: former results}

A primary issue was to explain the results of the biological experiments of \cite{filipski2002, filipski2004, filipski2005} by comparing the dominant eigenvalues in different settings. The original approach consisted in modelling perturbation of circadian rhythms through a loss of time dependency. A natural choice of relevant time-independent coefficients was the one given by the time average over a period, that is, obtained by replacing $d_i(t,x), K_{i\rightarrow i+1}(t,x)$ by $\langle d_i(x)\rangle, \langle K_{i\rightarrow i+1}(x)\rangle$ (where we use the notation $\langle f\rangle =T^{-1}\int_0^T f(s)ds$ for the time average). The new system possesses also a dominant eigenvalue that we now denote by $\lambda_P$ ($P$ for Perron, referring to the Perron-Frobenius theorem for positive linear operators). As the averaged system is supposed to represent perturbed rhythms, we expected to obtain generically $\lambda_F\leq \lambda_P$ (i.e., perturbed cells grow faster). The first result \cite{CMP}  was almost exactly the opposite: if the transition rates are time-independent, then one can prove that $\lambda_F\geq \lambda_P$. Biologically, this means that if a circadian control is exerted {\em only} on death rates, then one expects that perturbations of this control will actually lower the cell population growth instead of enhancing it!  To understand better the next results it is useful to consider even more general models:
\begin{equation}\label{eq:renew_system}
\left\lbrace\begin{array}{l}
\partial_t n_i(t,x)+\partial_x n_i(t,x)+d_i(t,x) n_i(t,x)=0,\\[0.3cm]
n_i(t,0)=\sum_j\int_0^\infty B_{j\rightarrow i}(t,x)n_j(t,x)dx\end{array}\right.
\end{equation}
The previous systems of PDEs  (\ref{eq:cycle_cellulaire},\ref{eq:division}) are only particular cases of such a system. What is important here is that we just separated birth rates  $B_{j\rightarrow i}$ and death rates $d_i$. It is worth noticing that in models (\ref{eq:cycle_cellulaire},\ref{eq:division}) transition rates contribute to both rates! The result of \cite{CMP} can also be generalized as follows: if the birth rates are time-independent then 
$\lambda_F\geq \lambda_P$. This result was improved in \cite{CGP} in the following way: consider system (\ref{eq:renew_system}) where we replace death rates by their time average and birth rates by their \textbf{geometrical} time average $\exp(T^{-1}\int_0^T \log  B_{j\rightarrow i}(t,x) dt)$, if we denote $\lambda_g$ the dominant eigenvalue of this system and as before $\lambda_F$ is the dominant eigenvalue associated to the time-dependent system, then, again $\lambda_F\geq \lambda_g$. We let the reader remark that if we consider this averaged version of the cell cycle system (\ref{eq:cycle_cellulaire}), the transition coefficients will be arithmetically averaged in the PDE but geometrically averaged in the boundary terms. 
We summarize the comparison results of \cite{CMP,CGP} in the following table:

\begin{table} 
\begin{tabular}{|c | c | c| c|}
		\hline 
			Birth rates & Death rates &  Dominant  & Inequalities\\
			 &  &  eigenvalue & \\ \hline
			 && Floquet &\\
			$B_{j\rightarrow i}(t,x) $& $d_i(t,x)$ & $\lambda_F$& \\[0.2cm] \hline
			&& geometric &\\
			$\exp(T^{-1}\int_0^T \log B_{j\rightarrow i}(t,x)dt)$&$T^{-1}\int_0^T d_i(t,x)dt$ & $\lambda_g $& $\lambda_g\leq \lambda_F$ (\cite{CGP})\\[0.2cm] \hline 
			&& Perron&\\
				$T^{-1}\int_0^T  B_{j\rightarrow i}(t,x)dt$&$T^{-1}\int_0^T d_i(t,x)dt$ & $\lambda_P $& $\lambda_g\leq \lambda_P$ (\cite{Lepoutre_MMNP})\\[0.2cm] \hline
		\end{tabular}
		\caption{Summary of the previous results on the comparison of eigenvalues.}
		\end{table}
		
		We refer the reader to \cite{Lepoutre_MMNP} where it is shown that the inequality $\lambda_g\leq \lambda_P$ is a straightforward consequence of the arithmetic geometric inequality $$\exp\left(T^{-1}\int_0^T \log B_{j\rightarrow i}(t,x)dt\right)\leq T^{-1}\int_0^T  B_{j\rightarrow i}(t,x)dt.$$ Note that it has also been theoretically established in \cite{Lepoutre_MMNP} that there is no generic inequality between $\lambda_F$ and $\lambda_P$.

\subsection{A unifying convexity result}
Our main result is the following convexity theorem.
\begin{Theorem} \label{theorem:convexityth}

The dominant eigenvalue associated to System~(\ref{eq:renew_system}) is convex with respect to the death rates $d_i$ and geometrically convex with respect to the birth rates $B_{j\rightarrow i}$. \\
In other words, consider  two sets of time $T-$periodic coefficients for System~(\ref{eq:renew_system}), namely 
$$(B_{j\rightarrow i}^1,d_i^1)_{1\leq i,j\leq I},\quad (B_{j\rightarrow i}^2,d_i^2)_{1\leq i,j\leq I}, $$
and denote $\lambda_F^1,\lambda_F^2$ the associated dominant eigenvalues.  Then for any $\theta\in [0,1]$, if we denote by $$d_i^\theta=\theta d_i^1+(1-\theta)d_i^2,\qquad B_{j\rightarrow i}^\theta=(B_{j\rightarrow i}^1)^\theta(B_{j\rightarrow i}^2)^{1-\theta}$$ the death rates and birth rates, respectively, and by $\lambda_F^\theta$ the dominant eigenvalue associated to coefficients 
$
d_i^\theta, B_{j\rightarrow i}^\theta
$, then we have  
\begin{equation}
\lambda_F^\theta\leq \theta \lambda_F^1+(1-\theta)\lambda_F^2.
\end{equation}

\end{Theorem}

 \begin{table}	\label{tab:convexity}
\begin{tabular}{|c | c | c| c|}
		\hline 
		
			Birth rates & Death rates &  Dominant  & Inequalities\\
			 &  &  eigenvalue & \\ \hline
			$B^1_{j\rightarrow i}$ & $d_i^1$ & $\lambda_F^1 $& \\ \hline
			$B^2_{j\rightarrow i}$ & $ d_i^2$ & $\lambda_F^2 $& \\ \hline
			&&&\\
			$(B^1_{j\rightarrow i})^{\theta} (B^2_{j\rightarrow i})^{1-\theta}$&$\theta d_i^1+(1-\theta)d_i^2$ & $\lambda_F^\theta $& $\lambda_F^\theta \leq \theta\lambda_F^1+(1-\theta)\lambda_F^2$\\ \hline
			
		\end{tabular}
		\caption{Summary of Theorem \ref{theorem:convexityth}.}
		\end{table}

Notice again that if we want to apply this to the cell cycle model, then we have to replace
 $K_{i\rightarrow i+1}$ by $\theta K_{i\rightarrow i+1}^1+(1-\theta)K_{i\rightarrow i+1}^2 $ in the PDE and by  $(K_{i\rightarrow i+1}^1)^{\theta} (K_{i\rightarrow i+1}^2)^{1-\theta} $ in the boundary terms. For the convenience of the reader, we summarize the theorem in Table 2
 (see Appendix and \cite{Lepoutre_PHD} for the proof).\\
 		
		Our convexity result generalizes the previous result of \cite{CGP} since the inequality $\lambda_g\leq\lambda_F$ can be recovered from Theorem~\ref{theorem:convexityth},  by using Jensen's inequality. 
	Considering different means (arithmetic and geometric) in the infinitesimal and integral term still lacks biological foundation. Especially, in the cell cycle model (\ref{eq:cycle_cellulaire}), since the geometrical mean is smaller than the arithmetical mean, this introduces in the infinitesimal term an artificial death term (the difference between the two means). Mathematically, it is easier to see where those different means come from by rewriting the simplest renewal equation (\ref{eq:renew1}) as a delay equation on $n(t,0)$ :
	$$
	n(t,0)=\int_0^{\infty} B(t,x)e^{-\int_0^x d(t-x+s,s)ds}n(t-x,0)dx.
	$$ 
	On this formulation, one can consider that only geometrical averages are used (on $B(t,x)e^{-\int_0^x d(t-x+s,s)ds}$), the geometrical averaging of $e^{-\int_0^x d(t-x+s,s)ds}$ is then equivalent to the arithmetical averaging of $d$.  This gives at least a mathematical, if not biological, justification for the introduction of a geometric, rather than arithmetic,  average for the birth rate.
		\subsection{An argument in favor of the indirect influence of the circadian control on tumour cells}
		
		An apparent default of this theoretical construction is that we do not explain in a direct way (as we would have expected)  the results of Filipski's experiments on tumor growth enhancement by circadian disruption, but as an interesting gain it provides a theoretical justification for the concept of chronotherapeutics. More precisely:
		
		An obvious limit of this approach is that it finally does not explain the experimental results of \cite{filipski2002, filipski2004, filipski2005}. Averaging the coefficients seems not to be the best way  to represent perturbations of circadian rhythms and we propose below alternative tracks to explain them.
		
		 However, as far as chronotherapeutics is concerned,  these comparisons are completely in keeping with the idea that the major effect of chronotherapy consists in minimizing toxicity on healthy tissues. If one considers that tumor tissues are less sensitive to circadian rhythms (in particular because one of the hallmarks of cancer cells is that they are insensitive to antigrowth factors, according to \cite{hanahan-weinberg}; this insensitivity to antigrowth factors is likely to extend to circadian inputs), we can focus only on circadian effects of drugs on healthy tissues. The goal of a therapy may then be thought of as to maximize the growth of healthy tissues that fight against a tumor, that is, to maximize $\lambda_F$, for healthy tissues only, of course.
		 
		  A first implication of the convexity inequality can be obtained if we model a therapy by an effect on death rate. Consider that the death rates $d_i(t,x)$ correspond to the periodic delivery of a fixed dose of a given drug at any time $0, T,2T,\dots$. We assume that  a phase-shifted delivery schedule at times $\varphi,T+\varphi, \dots $ will modify the death rates to $d_i(t+\varphi,x)$ and we denote the corresponding dominant eigenvalue by $\lambda(\varphi)$. A possible approach would consist in delivering the dose of drug uniformly over a period. It would lead to death rates $T^{-1}\int_0^t d_i(t,x)dt$. We denote the corresponding dominant eigenvalue by $\lambda_u$. Notice that neither the transition rates nor the birth rates are affected here by the drug (which is of course a strong assumption). The convexity inequality tells us that 
		$$
		\lambda_u\leq T^{-1}\int_0^T \lambda(\varphi)d\varphi.
		$$
		Hence, the phase may always be chosen in such a way that the periodic treatment is less toxic than its equivalent constant treatment (i.e., with the same daily dose).
		
		To illustrate this convexity inequality, we show two extreme cases of the profiles of $\lambda_F$ as a function of the phase shift: low advantage of the periodic treatment, but for a long range of phases, or else high advantage, but for a short range. Comparing the phase ranges in these two extreme cases with the equiprobabilistic (of being beneficial or detrimental) measure $T/2$ in the choice of the phase, one obtains the two following schematic situations:
		\begin{itemize}
		\item[-] the periodic delivery is for most of the time schedules (indexed by the phase shift parameter $\varphi$) more efficient than the constant treatment (mathematically $|\{\varphi, \; \lambda(\varphi)\geq\lambda_u\}|\geq T/2$, i.e., more than half of the available periodic delivery shifts do better than constant delivery); the gain is potentially low, but then highly probable even if the phase is chosen at random
		\item[-] the advantage obtained when $\lambda(\varphi)\geq \lambda_u$ is high for some $\varphi$: the gain is high, but the therapeutic window is narrow.
		\end{itemize}
		   
\begin{figure}
	\centering
		\includegraphics[width=0.75\textwidth]{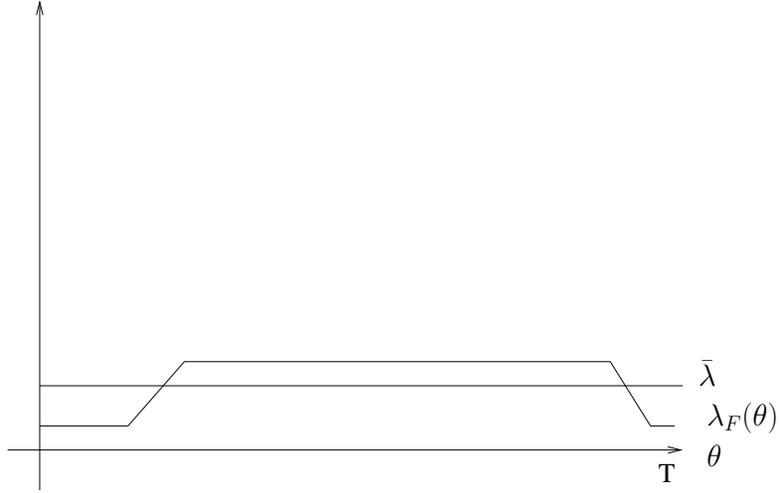}
	\caption{A schematic case $|\{\theta,\quad  \lambda(\theta)\geq \lambda_u\}| \geq T/2$. Here $\lambda(\theta)=\lambda_F(\theta)$ is the Floquet eigenvalue, as a function of phase shift $\theta$, and $\lambda_u=\overline{\lambda}$ is the Perron eigenvalue. }
	\label{fig:config_base_grande}
\end{figure}

\begin{figure}
	\centering
		\includegraphics[width=0.75\textwidth]{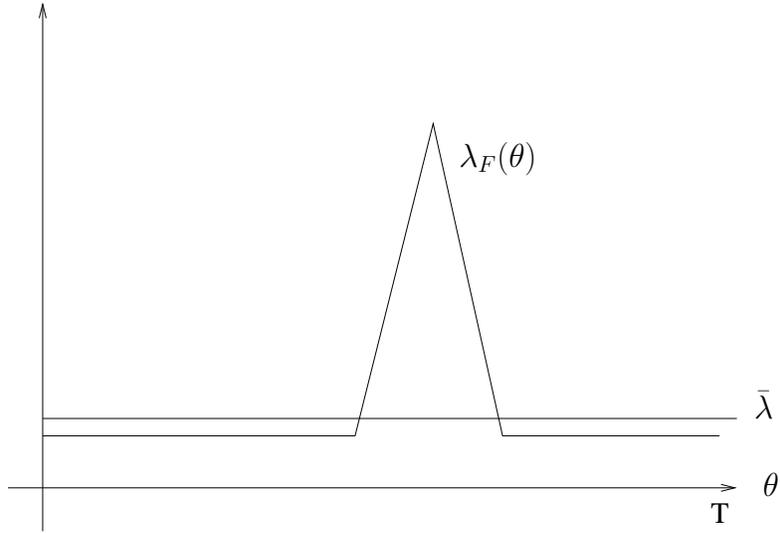}
	\caption{A schematic case $|\{\theta,\quad  \lambda(\theta)\geq \lambda_u\}| < T/2$, in which a significant advantage over constant infusion is obtained by chronotherapy. Same notations as in Figure \ref{fig:config_base_grande}.}
	\label{fig:config_base_petite}
\end{figure}

		These configurations are represented in Figure~\ref{fig:config_base_grande} and in Figure~\ref{fig:config_base_petite} (these figures being here only meant as illustrative sketches, and not results of simulations). Last but not least, this result also implies that if a periodic delivery may be more toxic (if there exists $\varphi$ such that $\lambda(\varphi)<\lambda_u$) then it also necessarily may be less toxic (meaning that there exists $\varphi$ such that $\lambda(\varphi)>\lambda_u$). 
		
		The therapeutic window may be (in some molecular way that remains to be elicited) genetically determined for a given individual and a given drug. This might explain puzzling situations that are found in the clinic, where it has been observed, for instance, that chronotherapy of colorectal cancer delivered at a commonly chosen peak phase resulted, by comparison with constant delivery, in significantly higher survival in males, but in significantly lower survival in females! (\cite{Sylvie, Sylvie2}). Our result does not explain such observation; it suggests however  that there should always be a therapeutic window that improves survival. Therefore, lower survival in women could be the result  of a presently unadapted chronotherapeutic protocol in females (so far, to our best knowledge, the same for men and for women). Our result suggests that, as in males, there should exist a best chronotherapeutic schedule for females, leading to a significant advantage in survival.

\begin{figure}
	\centering
		\includegraphics[width=0.75\textwidth]{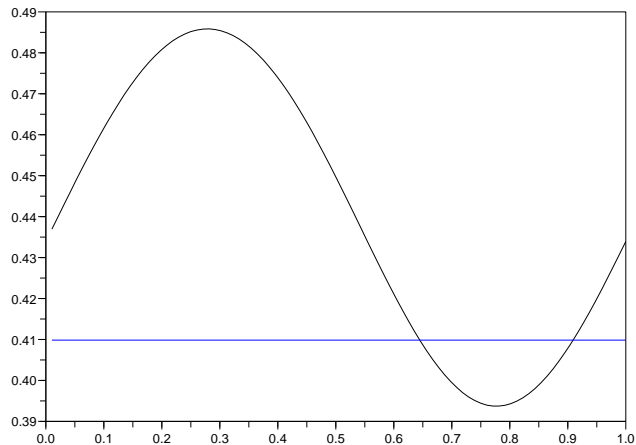}
	\caption{Possible advantages of chronotherapy: simulation for the coefficients given in Table 3 and comparison between first eigenvalues $\lambda(\varphi)$ (sinewave-like)  as a function of phase shift $\varphi$ in abscissae and $\lambda_u$ (constant line), for healthy tissue growth.}
	\label{fig:convexite1}
\end{figure}

\newpage
 We provide also a numerical illustration for this result. We consider a three-phase cell cycle model having the following structure:
$$ 
\left\lbrace
 \begin{array}{l}
 \partial_t n_i(t,x)+\partial_xn_i(t,x)+\left[d_i(t)+\psi_i(t)\chi_{[a_i,\infty[}(x)\right]n_i(t,x)=0,\quad 1\leq i\leq 3\nonumber, \\[0.3cm]\nonumber
   n_{i+1}(t,0)=\psi_i(t)\int_{a_i}^\infty n_i(t,x)dx,\quad i\leq 2,\nonumber\\[0.3cm]
 n_1(t,0)=2\psi_3(t)\int_{a_3}^\infty n_3(t,x)dx.\nonumber
 \end{array}\right.$$
 That is, we take system \ref{eq:cycle_cellulaire} with $I=3$, $d_i(t,x)=d_i(t)$ and $K_{i\rightarrow i+1}(t,x)=\psi_i(t)\chi_{[a_i,\infty[}(x)$, where $\chi_I$ stands for the indicator (characteristic) function of interval $I$. In a numerical experiment presented in Figure \ref{fig:convexite1}, the period being $1$, we have chosen a set of coefficients that is summarized in Table 3.
 \begin{table}
 \begin{tabular}{|c|c|}
 \hline
 $a_1$ & $10/24$\\
 $a_2 $& $10/24$\\
 $a_3$ & $2/24$\\
 \hline
 \end{tabular}\hspace{1cm}
  \begin{tabular}{|c|c|}\hline
 $\psi_1 $&$10*(1+0.8\cos(2\pi t))$\\
 $\psi_2 $& $10*(1-0.8\cos(2\pi t))$\\
 $\psi_3 $& 10\\
 \hline
  \end{tabular}\hspace{1cm}
  \begin{tabular}{|c|c|}\hline
 $d_1$ & $0$\\
  $d_2 $& $\cos^6(\pi t)$\\
  $d_3$ & $0$\\
  \hline
\end{tabular}
\caption{Coefficients in the numerical illustration of Figure \ref{fig:convexite1}}
\end{table}

\section{Possible alternative explanations for the initial observation}
\label{Alternative}

We have already mentioned that a possible explanation of the results of \cite{filipski2002, filipski2004, filipski2005} is that the main effect of a circadian control is exerted on healthy cells rather that on cancer cells, that are supposed to be less responsive to regulation factors, and from a chronotherapeutic point of view, that less unwanted toxicity on healthy cells results from a periodic, rather than constant, drug delivery schedule.

Another way to imagine a periodic control on the cell division cycle is to think (following an idea first developed in \cite{romond}) that it may be exerted on transition rates in a differentiated way, i.e., for instance with a common circadian input resulting in phase-opposed gate opening at the $G_1/S$ and $G_2/M$ transitions. This idea also tries to follow the biological observation that circadian gene expression of  Bmal1, that controls $G_2/M$ gate opening through Wee1 and Cyclin B-Cdk1, is in antiphase with the expression of another main circadian gene, Per2,  that controls $G_1/S$ gate opening through p21 and Cyclin E-Cdk2. In an unpublished simulation study on a 3-phase model \cite{Emilio1phase}, we have obtained that such phase-opposed periodic control between $G_1/S$ and $G_2/M$ transitions physiologically results in lower eigenvalues by comparison with constant transition rates, i.e., no control at all. That is, we always obtained the result $\lambda_F\leq\lambda_P$, contrarily to the main convexity result exposed above. By `physiologically', we mean here that we assumed a constant mitosis time, $1/24$ of the total cell cycle time, which is generally admitted by biologists of the cell cycle. In the same simulations, we tested completely freely varying times for the other two phases, $G_1$ and $S-G_2$, with total cell cycle duration time maintained constant.  Note anyway that an important difference in this alternative setting with the above framework is that we used only arithmetic, never geometric averaging, for the transition coefficients of the uncontrolled system. How can this discrepancy be explained (provided that one admits the validity of these simulation results, that are not accounted for by any theorem so far)? Independently of the choice for the averaged system of the arithmetic or of the geometric mean, for which we lack definitive biological justifications, this alternative explanation might be related to a physiological, as opposed to pathological, situation. In this  situation, circadian controls on healthy tissues are normally exerted without being hampered. The experiments of Filipski et al. \cite{filipski2002, filipski2004, filipski2005} occur on the contrary in a very pathological situation describing a fast growing tumor and its surrounding healthy tissue, that may be considerably perturbed, e.g., by tumor-emitted cytokines. This is only speculation so far, and more biological experiments remain to be done to confirm or infirm this hypothesis. In particular, the relationships that may exist between tissue synchronization with respect to cell cycle phases, in healthy and in tumour tissue, and phases of essential circadian proteins as Bmal1 and Per2 may be critical and should be investigated to that purpose.

\section{Conclusion}
\label{conclusion}

We use a general framework of physiologically structured PDEs to describe the cell division cycle in proliferating cell populations and its control by periodic inputs. Our motivation comes in particular from the knowledge of inputs from physiological circadian clocks on the cell divisin cycle, but also from experimental results of cancer chronotherapeutics  (periodic drug delivery). We have also given a possible theoretical justification for the success of cancer chronotherapeutics. Our mathematical results lead us to propose a simple, but not immediately patent, explanation  to account for the initial challenging biological observation (enhancement of tumor growth by circadian clock disruption), provided that we admit that it is not the tumor, but rather the healthy tissue that fights against it, that is the object of a perturbed circadian control. Nevertheless, this speculation still needs to be supported both by further experiments and by more elaborate, physiologically based, mathematical models.

\appendix

\section{Proof of the convexity result}
 \label{sec:convexity}
We restrict the proof of the main convexity result to the simple case of a single renewal equation. This result can be understood as a generalization of a famous result of Kingman \cite{Kingman} on nonnegative matrices. This result can be summarized as follows : if we denote $\rho$ the Perron root of nonnegative matrices, then $\log \rho$ is convex with respect to diagonal coefficients and geometrically convex with respect to the offdiagonal coefficients. This explains the different treatment for birth and death coefficients : death which impact is local is treated as a diagonal coefficient whereas the birth rate has a nonlocal impact and therefore is treated as an offdiagonal coefficient. No new difficulties arise in the generalization to systems. We study the behavior of the renewal equation, which is of course a particular case of (\ref{eq:renew_system}):
\begin{equation}\label{eq:renew1}
\left\lbrace\begin{array}{l}
\partial_t n(t,x)+\partial_x n(t,x) +d(t,x)n(t,x)=0,\\[0.3cm]
n(t,0)=\int_0^\infty B(t,x)n(t,x)dx.
\end{array}\right.
\end{equation}
where the birth rate $B$ and the death rate $d$ are taken nonnegative and $T-$periodic with respect to time $t$. The growth rate can be defined as the unique real $\lambda_F$ such that there exists two nonnegative  associated eigenfunctions $(N,\phi)$ satisfying:
\begin{equation}\label{eq:renew1_eigen}
\left\lbrace\begin{array}{l}
\partial_t N(t,x)+\partial_x N(t,x) +\left[d(t,x)+\lambda_F\right]N(t,x)=0,\\[0.3cm]
N(t,0)=\int_0^\infty B(t,x)N(t,x)dx,\\[0.3cm]
N(t+T,x)=N(t,x),\quad N\geq 0, N\not=0,\\[0.3cm]
-\partial_t \phi(t,x)-\partial_x \phi(t,x) +\left[d(t,x)+\lambda_F\right]\phi(t,x)=B(t,x)\phi(t,0),\\[0.3cm]
\phi(t+T,x)=\phi(t,x),\quad\phi>0.
\end{array}\right.
\end{equation}
If such eigenelements do not exist (which can happen if the birth coefficient $B(t,x)$ vanishes for some values oft and x ), then we can define $\lambda_F$ as the infimum of real $\mu$ such that there exist a positive dual subeigenfunction $\phi_\mu$ satisfying:
\begin{equation}\label{eq:renew1_eigen_sub}
\left\lbrace\begin{array}{l}
-\partial_t \phi_\mu(t,x)-\partial_x \phi_\mu(t,x) +\left[d(t,x)+\mu\right]\phi_\mu(t,x)\geq B(t,x)\phi_\mu(t,0),\\[0.3cm]
\phi_\mu(t+T,x)=\phi_\mu(t,x),\quad\phi_\mu>0.
\end{array}\right.
\end{equation}
This relaxed definition of the Floquet eigenvalue is inspired by the Collatz-Wielandt characterization of the dominant root arising in Perron-Frobenius theory (see \cite{HornJohnson} for instance).
Note that if $(\mu,\phi_\mu)$ satisfies (\ref{eq:renew1_eigen_sub}), then we have, 
$$
\frac{d}{dt}\int_0^\infty n(t,x)e^{-\mu t}\phi_\mu(t,x)dx\leq 0,
$$ 
therefore 
$$
\int_0^\infty n(t,x)\phi_\mu(t,x)dx\leq e^{\mu t}\int_0^\infty n(0,x)\phi_\mu(0,x)dx.
$$
In words, any solution grows slower than $e^{\mu t}$ in a weighted space. Note that when eigenelements exist the two notions coincide and the same computation leads to 
$$
\int_0^\infty n(t,x)\phi(t,x)dx=e^{\lambda_F t}\int_0^\infty n(0,x)\phi(0,x)dx.
$$
It justifies the idea that solutions grow like $e^{\lambda_F t}$. 
The proof of the theorem is mainly based on the following lemma
\begin{lemma}
Given two sets of $T-$ periodic coefficients $B_1,B_2,d_1,d_2$, if we can find  $(\mu_1, \phi_{\mu_1}^1)$ and $(\mu_2, \phi_{\mu_2}^2)$ satisfying (\ref{eq:renew1_eigen_sub}) (for $d=d_1,B=B_1$ and $d=d_2,B=B_2$ respectively) ,  then we have for any $\theta\in[0,1]$, denoting $\phi^\theta=(\phi_{\mu_1}^1)^\theta (\phi_{\mu_2}^2)^{1-\theta}$, $B^\theta=(B_1)^\theta (B_2)^{1-\theta}$ and 
$d^\theta=\theta d_1+(1-\theta)d_2$,
 \begin{equation}
\left\lbrace\begin{array}{l}
-\partial_t \phi^\theta(t,x)-\partial_x \phi^\theta (t,x) +[\theta\mu_1+(1-\theta)\mu_2 +d^\theta ]\phi^\theta(t,x)\geq B^\theta(t,x)\phi^\theta(t,0),\\[0.3cm]
\phi^\theta (t+T,x)=\phi^\theta(t,x),\quad\phi^\theta>0.
\end{array}\right.
\end{equation}
\end{lemma}
\textbf{Proof.} As $\phi^i_{\mu_i}$ is positive, we can write the equation on $\log\phi^i_{\mu_i}$ (we just have to divide by $\phi^i_{\mu_i}$):

\begin{eqnarray}
\!\!\!\!\!\!\!\!\!\!\!\!\!\!-\partial_t \log\phi^1_{\mu_1}(t,x)-\partial_x \log\phi^1_{\mu_1}(t,x) +(d_1(t,x)+\mu_1)&\!\!\geq&\!\! B_1(t,x)\frac{\phi^1_{\mu_1}(t,0)}{\phi^1_{\mu_1}(t,x)},\label{1}\\
\!\!\!\!\!\!\!\!\!\!\!\!\!\!-\partial_t \log\phi^2_{\mu_2}(t,x)-\partial_x \log\phi^2_{\mu_2}(t,x) +(d_2(t,x)+\mu_2)&\!\!\geq&\!\! B_2(t,x)\frac{\phi^2_{\mu_2}(t,0)}{\phi^2_{\mu_2}(t,x)}\label{2}
\end{eqnarray}
Thanks to the arithmetic geometric inequality, we have
$$
\theta B_1(t,x)\frac{\phi^1_{\mu_1}(t,0)}{\phi^1_{\mu_1}(t,x)}+(1-\theta)B_2(t,x)\frac{\phi^2_{\mu_2}(t,0)}{\phi^2_{\mu_2}(t,x)}\geq B^\theta(t,x)\frac{\phi^\theta(t,0)}{\phi^\theta(t,x)}.
$$
We also know that 
$$
\theta \log\phi^1_{\mu_1}(t,x)+(1-\theta) \log\phi^2_{\mu_2}(t,x)=\log \phi^\theta(t,x), \quad\theta d_1(t,x)+(1-\theta) d_2(t,x)=d^\theta(t,x).
$$
Noticing that, summing $\theta$(\ref{1})+$(1-\theta)$(\ref{2}) gives
$$
-\partial_t \log\phi^\theta(t,x)-\partial_x \log\phi^\theta(t,x) +(d^\theta(t,x)+\theta\mu_1+(1-\theta)\mu_2)\geq B^\theta(t,x)\frac{\phi^\theta(t,0)}{\phi^\theta(t,x)}.
$$
Multiplying by $\phi^\theta$, we get 
$$
-\partial_t \phi^\theta(t,x)-\partial_x \phi^\theta(t,x) +\bigg(d^\theta(t,x)+\theta\mu_1+(1-\theta)\mu_2\bigg)\phi^\theta(t,x)\geq B^\theta(t,x)\phi^\theta(t,0).
$$

\begin{cor}
With the above notations we have 
$$
\lambda_F^\theta\leq \theta\lambda_F^1+(1-\theta)\lambda_F^2.
$$
\end{cor}
\textbf{Proof.}
A first consequence of the proof of the previous lemma is that $$\lambda_F^\theta\leq \theta\mu_1+(1-\theta)\mu_2,$$ and from the definition we can choose $\mu_i\rightarrow \lambda_F^i$ and conclude.

\newpage
\bibliographystyle{alpha}
\bibliography{CGL4MCM2009_revising}







\end{document}